\documentclass[a4,11pt,reqno]{amsart}
\usepackage[english]{babel}    \usepackage[latin1]{inputenc}   \usepackage[T1]{fontenc}     \usepackage[french]{minitoc}
\usepackage[nice]{nicefrac}    \usepackage{latexsym,amsfonts}  %\usepackage[french]{babel}  %\usepackage{garamond}
\usepackage{graphics}    \usepackage{ulem}       \usepackage{hhline}    \usepackage{dsfont}    \usepackage{mathrsfs}
\usepackage{fancyhdr}    \usepackage{amsmath}    \usepackage{amssymb}   \usepackage{rotating}  \usepackage{fancybox}
\usepackage{color}       \usepackage{colortbl}   \usepackage{setspace}  \usepackage{enumerate} \usepackage{amsthm}
\usepackage{multicol}    \usepackage{pifont}     \usepackage{amsthm}    \usepackage{varioref}  \usepackage{textcomp}
\usepackage{lmodern}     \usepackage{mathpazo}   \usepackage{euscript}  \usepackage[pdftex]{hyperref}
\numberwithin{equation}{section}  \makeatletter\@addtoreset{equation}{section}
   \DeclareMathSymbol{\subsetneqq}{\mathbin}{AMSb}{36}
\newtheorem {theorem}{Theorem}[section]

\newtheorem {definition}[theorem]{Definition}

\newtheorem {proposition}[theorem]{Proposition}
\newcommand{\C}{\mathbb C}
\newcommand{\R}{\mathbb R}

\newcommand{\aF}{{{_2\digamma_1}}}
\newcommand{\bF}{{{_3\digamma_2}}}
\newcommand{\chp}{\nu}
\newcommand{\B}{\mathbb B^n}
\newcommand{\Berm}{\mathscr{B}_m}
\newcommand{\Berezin}{\mathscr{B}_0}
\newcommand{\gBerg}{\mathcal{A}_{m}^{2,\chp}(\B)}
\newcommand{\Berg}{{{\mathcal{A}}_{0}^{2,\nu}(\mathbb{B}^n)}}

         \newcommand{\fin}{\hfill $\square$}%\rule
            
\newcommand{\norm}[1]{\left\Vert#1\right\Vert} \newcommand{\set}[1]{\left\{#1\right\}} 
  \newcommand{\scal}[1]{\left<#1\right>}  
\begin{document}
%================================================================================================
%================================================================================================

\title[Magnetic Berezin transforms]{{A Formula representing magnetic Berezin transforms on the
Bergman ball of $\C^n$ as functions of the Laplace-Beltrami operator}}
\author{Allal GHANMI and Zouhair MOUAYN }
\address{{\bf (A.G.)} \quad  Department of Mathematics,  Faculty of  Sciences, P.O. Box 1014,
   Mohammed V University,  Agdal,  10 000 Rabat - Morocco  }
   %\ead{allalghanmi@gmail.com}
 %   \ead{intissar@fsr.ac.ma}
  %  \author{Z. Mouayn}
\address{{\bf (Z.M.)} \quad Department of Mathematics, Faculty of Sciences and Technics (M'Ghila), P.O.
Box 523, Sultan Moulay Slimane University, 23 000 B\'{e}ni Mellal - Morocco}
%\ead{mouayn@fstbm.ac.ma}

\begin{abstract}
We give a formula that represents magnetic Berezin transforms associated with
generalized Bergman spaces as functions of the Laplace-Beltrami operator on the
Bergman ball of $\C^{n}$. In particular, we recover the result obtained by J. Peeter [J. Oper. Theory, 24, 1990].

\quad


\bigskip

\bigskip
\end{abstract}

\maketitle
\section{\textsf{Introduction}}

The Berezin transform was introduced and studied by Berezin \cite{Berezin75}
for certain classical symmetric domains in $\C^n$
and next extensively considered in the context of Bergman, Hardy and
Bargmann-Fock spaces \cite{Berger,Englis,HKZ}.
It is closely related to Topelitz operators \cite{Zhu} and is useful in the quantization theory and corresponding $\star$-product on a suitable algebras of functions as well as in the correspondence principle  \cite{Berezin75,Berger}.

This transform is defined as follows. Let $\Omega$ be a domain of $\C^n$
with a Borel measure $\mu$ and $\mathcal{H}$ a closed subspace of $L^2(\Omega,d\mu)$ consisting of continuous
functions and  possessing a reproducing kernel $K(\cdot,\cdot)$. Then, the Berezin symbol $\sigma(A) $ of a
bounded linear operator ${A}$ on $\mathcal{H}$ is the function on $\Omega$ given by
$\sigma({A})(x) :=\scal{{A}e_{x},e_{x}}$, where $e_{x}(\cdot) :=K(\cdot,x) K(x,x)^{-1/2} \in\mathcal{ H}$.
Thus the Toeplitz operator $T_f$ with symbol $f\in L^{\infty}(\Omega)$ is the operator defined on $\mathcal{H}$ by
$T_f(\varphi) =P(f\varphi)$; $\varphi \in \mathcal{H}$,
where $P$ is the orthogonal projection from $L^2(\Omega,d\mu)$ into $\mathcal{H}$.
The Berezin transform associated to $\mathcal{H}$ is then defined to be the positive
self-adjoint operator $\mathscr{B}:=\sigma T,$ which turns out to be a bounded
operator on $L^2(\Omega, K(x,x)d\mu(x))$.

In this paper, we deal with the Bergman unit ball $\Omega=\B$ in $(\C^n,\scal{\cdot,\cdot})$ as domain endowed with
its Lebesgue measure $d\mu$. Since the Berezin transform can be defined
provided that there is a given closed subspace of $L^2(\B,d\mu)$ possessing a
reproducing kernel, we are concerned here with the  $L^2$-eigenspaces
\begin{equation}\label{gBerg}
\gBerg =\left\{ \varphi \in L^2(\B,(1-|z|^2)^{-n-1}d\mu) ; \quad H_{\chp}\varphi
=\epsilon^{\chp,n}_m\varphi \right\} ,
\end{equation}
associated to the discrete spectrum
\begin{equation}
\epsilon^{\chp,n}_m=4\chp(2m+n)-4m(m+n); \qquad m=0,1,2,\cdots, [\chp-\frac n2],
\end{equation}
of the Schr\"{o}dinger operator with uniform magnetic field on $\B$ given by
\begin{equation}\label{MagnSchr0}
H_{\chp}=-4(1-|z|^2) \left\{
\sum_{j=1}^{n}( \delta _{ij}-z_{i}\overline{z_{j}}) \frac{\partial^2}{\partial z_{i}\partial \overline{z_{j}}}
+\chp\sum_{j=1}^{n}\left( z_{j}\frac{\partial }{\partial z_{j}}-\overline{z_{j}}\frac{\partial }{\partial \overline{z_{j}}}\right)
+\chp^2\right\} + 4\chp^2,
\end{equation}
provided that $\chp >n/2$. Above $[x]$ denotes the greatest integer not exceeding $x$.
Note that for $m=0$,  the space $\Berg$ reduces further to be isomorphic to the weighted Bergman space of holomorphic functions $g$ on $\B$ satisfying the growth condition
\begin{eqnarray}\label{GCBerg}
 \int_{\B} |g(z)|^2(1-|z|^2)^{2\chp -n-1}d\mu(z) <+\infty .
\end{eqnarray}
The associated Berezin transform can be expressed as
\begin{equation}
\mathscr{B}[\varphi] (z) = \frac{(2\chp-n)\Gamma(2\chp )}{\pi^n\Gamma(2\chp -n+1)}
\int_{\B}\left(\frac{(1-|z|^2)(1-|\xi|^2)}{|1-<z,\xi>|^2}\right)^{2\chp } \varphi
(\xi) \left( 1-|\xi|^2\right) ^{-\left( n+1\right) }d\mu(\xi) .  \label{B0}
\end{equation}
Moreover,  it can be written as a function of the Laplace-Beltrami operator $\Delta_{\B}$
(\cite[p.182]{peeter}) as
%\begin{equation}\label{peetre}
% \Berezin  =\frac{\left| \Gamma(\alpha+\frac{n}{2}+1+i\Lambda) \right|^2}{\Gamma(\alpha+1) \Gamma(\alpha +n+1) } ,
%\end{equation}
\begin{equation}\label{peetre}
 \mathscr{B}  =\frac{\left| \Gamma\left(\alpha+1+\frac{n}{2} +\frac i2 \sqrt{-\Delta_{\B}-n^2}\right) \right|^2}{\Gamma(\alpha+1) \Gamma(\alpha +n+1) } ,
\end{equation}
where $\Gamma(x)$ being the usual Gamma function %, $\Lambda =\frac 12 \sqrt{-\Delta_{\B} -n^2}$
 and $\alpha =2\chp-n-1$ occurring in the measure weight of the Bergman space under consideration in \cite{peeter} by Peeter who derived \eqref{peetre} from a formula representing the Berezin transform on the Bergman ball as function of the $\Delta_{\B}$, which was first established by Berezin in \cite{Berezin2}.

% In general, formulae representing Berezin transform as function of the Laplace-Beltrami operator play a key role in the Berezin quantization \cite{Berezin75}.

Our aim is then to generalize \eqref{B0} and \eqref{peetre} to the case of the spaces in \eqref{gBerg}.
We precisely attach to each of these eigenspaces a Berezin transform $\Berm$ following the formalism described
above as
\begin{eqnarray}
\Berm [\varphi](z) &:&=\frac{\Gamma(n)m!(2(\chp-m)-n)\Gamma(2\chp-m)}{\pi^n \Gamma(n+m)\Gamma(2\chp-m-n+1)} \int_{\B}\left( \frac{(1-|z|^2)(1-|\xi|^2)}{|1-\scal{z,\xi}|^2}\right)
^{2(\chp-m) } \label{BerInt}\\
&&\qquad \times \left( P_{m}^{(n-1,2[\chp-m])}(1-2|\xi|^2) \right)^2  \varphi(\xi)( 1-|\xi|^2)^{-(n+1)}d\mu(\xi),
\nonumber
\end{eqnarray}
where $P_{m}^{(\alpha,\beta)}(\cdot)$ denotes the Jacobi polynomial \cite{GR}. Moreover, we have to
prove that such transform can also be expressed as a function of the Laplace-Beltrami operator
$\Delta_{\B}$ in terms of the $\bF$-sum as
%\begin{align}
%\Berm  &=\frac{(2(\chp-m)-n) \Gamma(n+m)}{m!\Gamma(2\chp -m-n+1) \Gamma(2\chp -n)}
%   \sum\limits_{j=0}^{2m}   (-2)^{j} \mathcal{B}(n+j,-\frac{n}{2}+i\Lambda+2(\chp-m)) A_{j} \label{***}
%\\ & \qquad \qquad \times {\bF }\left( \frac{n}{2}+i\Lambda,n+j,\frac{n}{2}+i\Lambda,\frac{n}{2}+i\Lambda+2(\chp-m) +j,n;1\right), \nonumber
%\end{align}
\begin{align}
\Berm  =C^{\chp,n}_m &\sum\limits_{j=0}^{2m}   (-2)^{j} A_{j}
    \frac{\Gamma(n+j)\Gamma\left(2(\chp-m)-\frac 12\left(n-i\sqrt{-\Delta_{\B} -n^2}\right) \right) }
    {\Gamma\left( 2(\chp-m) +j + \frac 12\left(n+i\sqrt{-\Delta_{\B} -n^2}\right)\right) } %\mathcal{B}\left( n+j,-\frac{n}{2}+\frac i2\sqrt{-\Delta_{\B} -n^2}+2(\chp-m) \right)
     \label{***}
\\ &  \times {\bF }\left[ \begin{array}{c}
                 \frac 12\left(n+i\sqrt{-\Delta_{\B} -n^2}\right), \, n+j, \, \frac 12\left(n+i\sqrt{-\Delta_{\B} -n^2}\right) \\
                 (\chp-m) +j+\frac 12\left(n+i\sqrt{-\Delta_{\B} -n^2}\right), \, n
               \end{array} \bigg| 1
\right],
%{\bF }\left( \frac 12\left(n+i\sqrt{-\Delta_{\B} -n^2}\right),n+j,\frac 12\left(n+i\sqrt{-\Delta_{\B} -n^2}\right);2(\chp-m) +j+\frac 12\left(n+i\sqrt{-\Delta_{\B} -n^2}\right),n;1\right),
\nonumber
\end{align}
where
\begin{align}
C^{\chp,n}_m= \frac{(2(\chp-m)-n) \Gamma(n+m)}{m!\Gamma(2\chp -n-m+1) \Gamma(2\chp -n)}
\end{align}
and
\begin{align}
A_{j}=2^{-j}\sum_{p=\max \left( 0,j-m\right) }^{\min \left( m,j\right) }%
%\binom{m}{j-p}\binom{m}{p}
\frac{(m!)^2\Gamma(2\chp-m) \Gamma(2\chp-m+j-p)}{(j-p)!(m+p-j)!p!(m-p)!\Gamma(n+j-p) \Gamma(n+p) }.
\end{align}
In particular, when $m=0,$ the transform \eqref{BerInt} reduces to the well
known Berezin transform  $\mathscr{B} $ in \eqref{B0} on the unit ball \cite{Berezin75}
and the expression \eqref{***} becomes the formula \eqref{peetre}.

The paper is organized as follows. In Section 2, we fix some notations and
we review briefly some needed tools
on the bi-holomorphic maps of the ball, the spectral function of the
Laplace-Beltrami operator and the Fourier-Jacobi transform of this operator. In Section 3, we are dealing
with the $L^2$-spectral theory of the Schr\"{o}dinger operator with uniform  magnetic
field on the unit ball and we are concerned in particular with  the reproducing kernels of
the eigenspaces associated with the discrete spectrum. In Section 4, we
construct for each of these eigenspaces a Berezin transform that we express
 in terms of the bi-holomorphic maps of the ball. In Section 5, we give a
formula representing the constructed Berezin transform as function of the
Laplace-Beltrami operator.

\section{Some notations and tools}
Here we fix some notations and we recall some needed tools we will be using such as the bi-holomorphic maps of the ball,
the spectral function of the Laplace-Beltrami operator and the Fourier-Jacobi transform.

\subsection{ The Moebius group of $\B$.}

Let $\C^n$ be endowed with its standard inner product
$
\scal{z,w} =z_{1}\overline{w_{1}}+z_{1}\overline{w_{2}}+\cdots +z_{1}\overline{w_{n}}
$
for every $z=(z_1,\cdots, z_n)$, $w=(w_1,\cdots,w_n)$ in $\C^n$, so that
$
|z| =\sqrt{\scal{z,z} }=\sqrt{|z_{1}|^2+ \cdots +| z_{n}|^2}.
$
We denote by $\B=\left\{ z\in \C^n:|z| <1\right\}$ the unit ball of $\C^n$. The Moebius group of $\B$, denoted $Aut (\B)$, i.e.
the set of the biholomorphic mappings  $\varphi: \B\longrightarrow \B$ (called also the
autormorphisms of $\B$) acts transitively on $\B$. For fixed $a\in \B$, let $P_{a}$ be the orthogonal projection of
$\C^n$ onto the subspace $<a>$ generated by $a$ and let $Q_{a}=I-P_{a}$ be
the projection onto the orthogonal complement of $<a>$. To be quite explicit, we have
$P_{a}(z)= {\scal{z,a} a}/{\scal{a,a} }$ if $a\neq 0$ and $P_{0}(z)=0$.
For every $a\in \B$, set $s=(1-|a|^2)^{1/2}$ and define
\begin{equation}
\varphi _{a}(z) =\frac{a-P_{a}(z) - (1-|a|^2)^{1/2}Q_{a}(z) }{1-\scal{z,a}}.
\end{equation}
 Then, the maps $\varphi _{a}$ have the following properties \cite{rudin}:
(i)  $\varphi _{a}\left( 0\right) =a$ \ \ and \ \ $\varphi _{a}(a) =0$, (ii)  $\varphi _{a}^{\prime}(0) =-s^2P_a-sQ_a$ \ \ and \ \ $\varphi_{a}^{\prime}(a) =-\frac{1}{s^2}P_a-\frac{1}{s}Q_a$, (iii) $\varphi _{a}$ is an involution: $\varphi _{a}\left(\varphi _{a}(z) \right) =z$,
 (iv)  $\varphi _{a}: \overline{\B}\rightarrow \overline{\B}$ is an homeomorphism  and belongs to $Aut(\B)$, (v) The Jacobian $\frac{|Dw|}{|D\xi|}$  of the transformation $\xi =\varphi_{a}(w)$ is equal to
\begin{equation}\label{Jacobian}
\frac{\left| Dw\right| }{\left| D\xi \right| }=\left( \frac{1-\left|
a\right|^2}{\left| 1-\scal{a,\xi } \right|^2}\right) ^{n+1}
\end{equation}
and (vi) The identity
\begin{equation}\label{Inva}
1-\scal{\varphi _{a}(w) ,\varphi _{a}(w')}
=\frac{\left( 1-\scal{a,a} \right) \left(1-\scal{w,w'} \right) }{\left( 1-\scal{w,a} \right) (1-\scal{a,w'} )}
\end{equation}
holds for every $w,w'\in \overline{\B}$.

The elements of $Aut(\B)$ are the maps $\varphi_{a}$; $a\in\B$, and their unitary transformations.
 Indeed, if $\psi \in Aut(\B) $ and $a=\psi ^{-1}(0)$, then there is a unique $U$ in the unitary group $U(n)$ such that $\psi =U.\varphi _{a}$.
 Moreover, in view of \eqref{Inva}, $\psi$ satisfies the identity
\begin{align}\label{IdentityInv}
1-\scal{\psi (w) ,\psi (w') } =%
\frac{\left( 1-\scal{a,a} \right) \left( 1-\scal{w,w'} \right) }{\left( 1-\scal{w,a} \right)
(1-\scal{a,w'} )}; \qquad w,w'\in \overline{\B}.
\end{align}
For our purpose, we will be concerned with the automorphism
\begin{align}
\varphi _{z}(w) :=\frac{A_{z}w+z}{1+\scal{w,z} }, \label{varphiz}
\end{align}
where $A_{z}$ is the $n\times n$ matrix given by
\begin{align}
A_{z}:=(1-|z|^2) ^{\frac{1}{2}}I_{n}+\frac{zz^{\ast }}{1+\sqrt{1-|z|^2}}, \label{varphiAz}
\end{align}
$I_n$ being the identity matrix.

\subsection{The Laplace-Belrami operator. }

Recall that the Bergman kernel of the unit ball is given by
 \begin{equation}
K(z,w) =\frac{n!}{\pi ^{n}\left( 1-\scal{z,w} \right) ^{n+1}},
\end{equation}
so that one defines
 \begin{equation}
g_{ij}=\partial _{i}\overline{\partial }_{j}\left(\log K(z,z)\right) =\frac{n+1%
}{\left( 1-|z|^2\right)^2}\left( \left( 1-|z|
^2\right) \delta _{ij}-\overline{z}_{i}z_{j}\right),
 \end{equation}
where we have used the notation $\partial _{k}$ and $\overline{\partial }_{k}$ to mean
\begin{equation}
\partial _{k}=\frac{\partial }{\partial z_{k}}=\frac{1}{2}\left( \frac{%
\partial }{\partial x_{k}}-i\frac{\partial }{\partial y_{k}}\right) \mbox{ and } \ \ \overline{\partial }_{k}=\frac{\partial }{\partial \overline{z}_{k}}=%
\frac{1}{2}\left( \frac{\partial }{\partial x_{k}}+i\frac{\partial }{%
\partial y_{k}}\right)
\end{equation}
with $z_{k}=x_{k}+iy_{k}\in \C$.
Therefore, the Bergman metric on $\B$ defined through
\begin{align}
ds_{z}^2 := \sum\limits_{i,j} \partial _{i}\overline{\partial }_{j}\log K\left(K(z,z)\right) \xi _{i}\overline{\xi }_{j}
 =(1-|z|^2) ^{-2}\sum\limits_{i,j=1}^{n}\left((1-|z|^2) \delta _{i,j}+\overline{z_{i}}z_{j}\right) dz_{i}\otimes d\overline{z_{j}},
\end{align}
so that the Bergman ball $(\B,ds^2)$ carries a K\"ahlerian structure.
The corresponding invariant volume measure  is
\begin{equation} \label{vm}
dg(z) =K(z,z) dV(z) =\frac{n!}{\pi ^{n}}\frac{dV}{\left( 1-|z|^2\right)^{n+1}}
=(1-|z|^2) ^{-(n+1)} d\mu(z)
\end{equation}
 and the geodesic distance $d(z,w)$ reads
\begin{equation}\label{Dist}
\cosh^2\left( d(z,w) \right) =\frac{|1-\scal{z,w}|^2}{(1-|z|^2)(1-|w|^2)}.
\end{equation}
 The Laplace-Beltrami operator associated to $(\B,ds^2)$ is defined to be
\begin{equation} \label{LaplBeltrami}
\Delta_{\B}=4\sum\limits_{i,j}g^{ij}\partial _{j}\overline{\partial }_{i}
=-4(1-|z|^2)\sum\limits_{i,j}\left( \delta _{ij}-z_i\overline{z}_{j}\right) \partial_{i}\overline{\partial }_{j};
\end{equation}
$(g^{ij})$ being the inverse matrix of $(g_{ij})$.
Note that $\Delta_{\B}$ is invariant under the action of the  group $SU(n,1)$ and
its spectrum $\sigma(\Delta_{\B})$ in the Hilbert space $L^2(\B,(1-|z|^2)^{-n-1}d\mu)$
is purely continuous and equals to $]-\infty, -n^2]$.
Instead of $\Delta_{\B}$, we have to deal with the non-negative elliptic self-adjoint operator
$-L_n:= -\Delta_{\B} -n^2$ so that $\sigma(-L_n)=[0,+\infty[$. Hence, for a given suitable real valued function $f: \R \longrightarrow \R$,
%is a measurable function,
 the operator $f\left( -L_{n}\right) $ is defined by
  \begin{align} \label{FctLn}
    f(-L_n)[\varphi](z) := \int_{\B} K_f(z,w) \varphi(w)(1-|w|^2)^{-n-1}d\mu(w),
    \end{align}
  whose the kernel function has the integral representation
\begin{align} \label{SpecKernel}
  K_f(z,w) = \int_{0}^{+\infty} \Psi(z,w;\lambda) f(\lambda) d\lambda,
\end{align}
 where the spectral kernel is given by (\cite[p.13]{Oueld}):
\begin{equation}
\Psi(z,w;\lambda) =   \frac{|\Gamma(\frac{n+i\lambda}2)|^4}{4\pi^{n+1}\Gamma(n)|\Gamma(i\lambda)|^2}
\, {\aF}\left[ \begin{array}{c} \frac{n+i\lambda }{2},\frac{n-i\lambda }{2}\\ n \end{array} \bigg | -\sinh^2 (d(z,w))\right] ,\label{BLBOueld}
\end{equation}
in terms of the ${_2\digamma_1}(\cdot)$-Gauss hypergeometric function \cite{GR}.
 %\begin{equation}
%{_2\digamma_1}(x)=\sum_{n=0}^\infty \frac{(a)_n(b)_n}{(c)_n} \frac{x^n}{n!}; \qquad |x|<1.
% \end{equation}

\subsection{The Fourier-Jacobi Transform. }

%In this section, we recall some known facts on the Fourier-Jacobi transform.
%For more information see Koornwinder's paper \cite{Koornwinder}  and the references therein.
Let us recall  that the Jacobi functions are defined by
\begin{align}\label{JacFun}
\phi _{\lambda }^{(\alpha,\beta) }(t)
:=\aF \left[ \begin{array}{c} \frac{1}{2}\left( \alpha +\beta +1+i\lambda
\right) ,\frac{1}{2}\left( \alpha +\beta +1-i\lambda \right) \\ \alpha +1 \end{array} \bigg | -\sinh^2(|t|)\right] .
\end{align}
The particular case of $t=i\theta$ and $\lambda=i(2n+\alpha+\beta+1)$ leads to
\begin{eqnarray}
\phi _{i\left( 2n+\alpha +\beta +1\right) }^{(\alpha,\beta)
}\left( i\theta \right) &:=&\aF \left[ \begin{array}{c} -n,n+\alpha +\beta
+1-i\lambda \\ \alpha +1 \end{array} \bigg |\sin^2(|\theta|) \right] \\
&=&\frac{n!}{\left( \alpha +1\right) _{n}}P_{n}^{(\alpha,\beta) }\left( \cos (|2\theta|) \right),
\end{eqnarray}
 a normalized Jacobi polynomial. We have the following special cases
%\begin{align}
%& \phi _{\lambda }^{\left( -\frac{1}{2},\frac{1}{2}\right) }(t)
%=\cos (\lambda t) \label{Case1}\\
%&\phi _{\lambda }^{\left( 0,0\right) }(t) =P_{\frac{1}{2}\left(
%i\lambda -1\right) }\left( \cosh (2t)\right), \label{Case2}
%\end{align}
$\phi _{\lambda }^{\left( -\frac{1}{2},\frac{1}{2}\right) }(t)
=\cos (|\lambda t|)$ and $\phi _{\lambda }^{\left( 0,0\right) }(t) =P_{\frac{1}{2}\left(
i\lambda -1\right) }\left( \cosh (|2t|)\right)$,
where $P_{\nu }(\cdot)$ is the Legendre function.
%We have also the following limit relationship to the Bessel function $J_{\alpha }(\cdot) $:
%\begin{align}
%\lim_{r\rightarrow \infty }\phi _{r\lambda }^{(\alpha,\beta)
%}\left( r^{-1}t\right) =2^{\alpha }\Gamma(\alpha+1) \left(
%\lambda t\right) ^{-\alpha }J_{\alpha }\left( \lambda t\right) . \label{Case3}
%\end{align}
Now, let us assume that $\alpha >-1$ and $\beta \in \Bbb{R\cup }i\R$, and set
\begin{align}\label{Deltaab}
\Delta _{\alpha ,\beta }(t) :=\left( 2\sinh (|t|)\right) ^{2\alpha+1}\left( 2\cosh(|t|)\right) ^{2\beta +1}.
\end{align}
Then the Jacobi function $\phi _{\lambda }^{(\alpha,\beta) }$ in \eqref{JacFun} is an even $%
C^{\infty }$-function on $\R$ satisfying the second order differential equation
 \begin{equation}
\left( \frac{d^2}{dt^2}+\frac{\Delta _{\alpha ,\beta }^{\prime }(t) }{\Delta _{\alpha ,\beta }(t) }\frac{d}{dt}+\lambda
^2+(\alpha+\beta+1)^2\right) \phi _{\lambda }^{(\alpha,\beta) }(t) =0
 \end{equation}
with $\phi _{\lambda }^{(\alpha,\beta)}(0) =1$, and is uniquely defined
by these properties.
Moreover, for $\alpha >-1$ and $\beta \in \left[ -\alpha -1,\alpha +1\right] \cup i%
\R$, we have the integral transform pair
\begin{align}
g(\lambda) &=\int\limits_{0}^{+\infty }f(t) \phi_{\lambda }^{(\alpha,\beta) }(t) \Delta _{\alpha,\beta }(t) dt \label{JacTransf}\\
f(t) &=(2\pi) ^{-1}\int_{0}^{+\infty }g(\lambda) \phi _{\lambda }^{(\alpha,\beta) }(t) \left| c_{\alpha ,\beta }(\lambda) \right|^{-2}d\lambda , \label{InvJacTransf}
\end{align}
where %$c_{\alpha ,\beta }(\lambda)$ stands for
\begin{align}
c_{\alpha ,\beta }(\lambda) :=\frac{2^{\alpha +\beta
+1-i\lambda }\Gamma(\alpha+1) \Gamma(i\lambda) }{%
\Gamma \left( \frac{1}{2}\left( \alpha +\beta +1+i\lambda \right) \right)
\Gamma \left( \frac{1}{2}\left( \alpha -\beta +1+i\lambda \right) \right) }. \label{ChandraT}
\end{align}
This establishes a one-to-one correspondence between the space of even $C^{\infty
}$-functions on $\R$ with compact support and its image under the
Fourier-cosine transform, as characterized by the classical Paley-Wiener
theorem. The mapping $f\mapsto g$, which is called the Fourier-Jacobi transform, can be extended
to an isometry of the Hilbert spaces between $L^2\left( \R%
_{+},\Delta _{\alpha ,\beta }(t) dt\right) $ and $L^2\left(
\R_{+},(2\pi) ^{-1}\left| c_{\alpha ,\beta }(\lambda) \right| ^{-2}d\lambda \right)$. Note that in the case of $\alpha >-1$ and $%
\left| \beta \right| >\alpha +1$, \eqref{JacTransf} remains valid except
that we have to add to the right hand side of the second formula in \eqref{InvJacTransf} the term
 \begin{equation}
\sum_{\lambda \in D_{\alpha ,\beta }}g(\lambda) d_{\alpha
,\beta }(\lambda),
 \end{equation}
where $D_{\alpha ,\beta }$ is a finite subset of $i\R$ and the weights $%
d_{\alpha ,\beta }(\lambda) $ can be given explicitly.
%\begin{remark}
%In view of the special cases \eqref{Case1}, \eqref{Case2} and \eqref{Case3}, the Fourier-Jacobi transform becomes
%the Fourier-cosine transform for $\alpha =\beta =-\frac{1}{2}$, the
%Mehler-Fock transform for $\alpha =\beta =0$, and converges to the Hankel
%transform under suitable changes of $\lambda$ and $t$.
%\end{remark}
%
%In this section, we recall some known facts on the Fourier-Jacobi transform.
 For more information  on the Fourier-Jacobi transform see Koornwinder's paper \cite{Koornwinder}  and the references therein.

\section{The generalized Bergman spaces $\mathcal{A}_{m}^{\chp}\left( \Bbb{B}^{n}\right) $}

Associated to the K\"ahlerian manifold $(\B,ds^2)$, modeling the complex hyperbolic space
 of negative holomorphic sectional curvature, there is a canonical vector potential
\begin{eqnarray}
\theta(z)=- \sqrt{-1} (\partial -\bar \partial)\mbox{Log} (1-|z|^2).
\end{eqnarray}
Thus, the magnetic Schr\"odinger operator, describing a charged particle moving in $\B$
  under the action of the derived magnetic field $d\theta$, can be defined as
\begin{eqnarray}
H_{\chp}:=  ( d + \sqrt {-1}\chp  \theta  )^{*} ( d + \sqrt {-1}\chp \theta).
\end{eqnarray}
It appears then as the Laplace-Beltrami operator in \eqref{LaplBeltrami}
  perturbed by a first order differential operator. Its   expression
in the complex coordinates is given by (\cite{GhInJMP}):
\begin{eqnarray} \label{LanLikeHam}
H_{\chp} = -4\Big(1-|z|^2\Big)\set{\sum _{i,j=1}^n(\delta _{ij}- z_i{\bar z_j})\frac{\partial ^2}{\partial z_i
\partial {\bar z_j}} +\chp \sum _{j=1}^n\left(z_j\frac{\partial}
{\partial z_j}-{\bar z_j}\frac{\partial}{\partial{\bar z_j}}\right) +\chp^2} +4\chp^2 .
\end{eqnarray}
%This second order differential operator can be related to the one considered in \cite{zhang,ahern96,BouIn} and
Different aspects of the spectral analysis of $H_{\chp}$ have been studied by many authors under different names and notations (like Maass Laplacians, Landau Hamiltonians, ... ), see for example \cite{zhang,ahern96,BouIn} for $n\geq 2$ and \cite{elstrodt,patterson} for $n=1$.
 For instance, note that $H_{\chp}$ is an elliptic densely defined operator on the Hilbert space  $L^2(\B,(1-|z|^2)^{-n-1}d\mu)$  and admits a unique self-adjoint realization denoted also by $H_{\chp}$.
Moreover, it is a known fact that the discrete spectrum of $H_{\chp}$ on $L^2(\B,(1-|z|^2)^{-n-1}d\mu)$ is not trivial whenever $\chp >n/2$ and consists of a finite number of infinitely degenerate eigenvalues of the form
\begin{eqnarray} \label{LandauLevels}
\epsilon^{\chp,n}_m:= 4\chp(2m+n)-4m(m+n)
\end{eqnarray}
for varying integer $m$ such that $0\leq m< \chp- n/2$. The corresponding $L^2$-eigenfunctions, i.e., the elements of
 the space $\gBerg$ in \eqref{gBerg}
 %the performed  $L^2$-Hilbert space
%\begin{eqnarray}
%\gBerg =\set{\varphi \in L^2 (\B,(1-|z|^2)^{-n-1}d\mu); \quad  H_{\chp} \varphi = \epsilon^{\chp,n}_m F }\label{HilbertB}
%\end{eqnarray}
can be written explicitly in terms of the Jacobi polynomials $P^{(\alpha,\beta)}_j(x)$ and the complex spherical harmonics
 $h^{p,q}(z,\bar z)$ \cite{folland,rudin}.  Namely, from \cite{zhang,GhInJMP} we check that the functions
\begin{eqnarray}
\psi^{\chp,m}_{p,q}(z) %=(1-|z|^2)^{\chp-m}  {\aF} (q-m, 2\chp +p-m ; n+p+q ; |z|^2)h^{p,q}_k(z,\bar z),
 %= \frac{(m-q)!\Gamma(n+p+q)}{\Gamma(n+m+p)}(1-|z|^2)^{\chp-m} P^{(n+p+q-1,2[\chp-m]-n)}_{m-q}(1-2|z|^2) h^{p,q}_k(z,\bar z),
:= (1-|z|^2)^{\chp-m} P^{(n+p+q-1,2[\chp-m]-n)}_{m-q}(1-2|z|^2) h^{p,q}(z,\bar z),
\label{L2eigenfunctions}
\end{eqnarray}
  constitute an orthogonal basis of the space in \eqref{gBerg}, for varying integers $p$ and $q$ such that $p=0,1,2, \cdots$, and  $q=0,1,2, \cdots ,m$, whose the square norms in $L^2(\B,(1-|z|^2)^{-n-1}d\mu)$ are given by
    \begin{eqnarray}
\norm{\psi^{\chp,m}_{p,q}}^2%_{L^2(\B,d\mu)}
=\frac{\pi^n\Gamma(n+m+p)\Gamma(2\chp -n-m-q+1)}
{2n!(2[\chp-m]-n) (m-q)!\Gamma(2\chp -m+p)} . \label{L2normb}
\end{eqnarray}
Furthermore, the space $\gBerg$ is a reproducing kernel Hilbert space, i.e.,
for every fixed $z\in \B$, there exists a function $K^{\chp}_{m,z}$ belonging to this space
such that
 \begin{equation}
f(z) = \scal{f, K^{\chp}_{m,z} } = \int_{\B} K^{\chp,n}_{m}(z,w)f(w)(1-|z|^2)^{-n-1}d\mu(w)
 \end{equation}
 for every $f \in \gBerg$,  where we have set $K^{\chp,n}_{m}(z,w)=\overline{K^{\chp}_{m,z}(w)}$.
 Its expression is given in terms of the hyperbolic distance in \eqref{Dist} as
\begin{align}
K^{\chp,n}_{m}(z,w)=\gamma^{\chp,n}_m \left(\frac{1- \overline{\scal{z,w}}}
{1-\scal{z,w}}\right)^{\chp}&\left(\cosh(d(z,w))\right) ^{-2(\chp-m)}  \label{rkB2}\\ & \times   P_m^{(n-1,2[\chp-m]-n)}\left(1-2\tanh^2(d(z,w))\right),\nonumber
\end{align}
where %the constant $\gamma^{\chp,n}_m$ is given explicitly  by
\begin{eqnarray}
\gamma^{\chp,n}_m= \frac{(2(\chp-m) -n)\Gamma(2\chp -m)}{\pi^n\Gamma(2\chp -m-n+1)}.\label{crkb}
\end{eqnarray}
For $m=0$, the $L^2$-eigenspace associated to the bottom
eigenvalue $\epsilon^{\chp,n}_0=4\chp n$  is isomorphic to the usual
weighted Bergman Hilbert space of holomorphic functions $g$ on $\B$ such that \eqref{GCBerg}
%\begin{eqnarray*}
%  \int_{\B} |g(z)|^2(1-|z|^2)^{2\chp -n-1}d\mu(z) <+\infty ,
%\end{eqnarray*}
holds. The corresponding reproducing kernel reads
 \begin{eqnarray}
K^{\chp,n}_{0}(z,w)=\frac{(2\chp -n)\Gamma(2\chp )}{\pi^n\Gamma(2\chp -n+1)}
\Big(\cosh(d(z,w))\Big) ^{-2\chp }.
\end{eqnarray}

\begin{definition}
Assume that $\chp > n/2$ and let $m$ be an integer such that $0\leq m <\chp - n/2$. Then the space $\gBerg$ is
 called here a generalized Bergman space of index $m$.
\end{definition}

\section{Generalized Berezin transforms on the unit ball}
We adopt to the formalism presented in Section 1 to the case of the Bergman ball $\Omega=\B$ as domain of $\C^n$ and $\gBerg$ as a closed subspace  of $L^2(\B,(1-|z|^2)^{-n-1}d\mu )$. Then, the Berezin symbol of a bounded
operator  $A$ on $\gBerg$ is the function
 \begin{equation}
\sigma _{m}(A) (z) :=\scal{A\left( \widetilde{%
K_{m}}\right) _{z}, \left( \widetilde{K_{m}}\right) _{z}}
_{L^2(\B,(1-|z|^2)^{-n-1}d\mu )}; \qquad z \in \B.
 \end{equation}
Here $(\widetilde{K_{m}})_{z}$ denotes the normalized reproducing kernel with evaluation at $z$, that is
 \begin{equation}
(\widetilde{K_{m}})_{z}(w) =\frac{K^{\chp,n}_{m}(z,w) }{\sqrt{%
K^{\chp,n}_{m}(z,z) }},\qquad w\in \B,
 \end{equation}
where $K^{\chp,n}_{m}(z,w)$ is given in \eqref{rkB2}.
% by
%\begin{align}
%K^{\chp,n}_{m}(z,w)=\gamma^{\chp,n}_m \left(\frac{1- \overline{\scal{z,w}}}
%{1-\scal{z,w}}\right)^{\chp}&\left(\cosh(d(z,w))\right) ^{-2(\chp-m)}  \nonumber \\ & \times   P_m^{(n-1,2[\chp-m]-n)}\left(1-2\tanh^2(d(z,w))\right),\nonumber
%\end{align}
%with
%\begin{eqnarray}
%\gamma^{\chp,n}_m= \frac{(2(\chp-m) -n)\Gamma(2\chp -m)}{\pi^n\Gamma(2\chp -m-n+1)}.\label{crkb}
%\end{eqnarray}
 Therefore, for given
 complex-valued function $\phi$ such that $\phi \widetilde{K_{m}} \in L^2(\B,(1-|z|^2)^{-n-1}d\mu )$,
 the Berezin transform of $\phi$ is defined to be the Berezin symbol of the Toeplitz operator $T_{m}(\phi ) $ with symbol $\phi $
 defined on $\gBerg $ by $T_{m}(\phi )[f]:=P_{m}[\phi f]$. That is,
\begin{align}
\Berm\left[ \phi \right] (z)  =\sigma _{m}(T_{m}(\phi))(z)
 =\scal{T_{m}(\phi )\left( \widetilde{K_{m}}\right) _{z},
\left( \widetilde{K_{m}}\right) _{z}} =\scal{P_{m}\left( \phi .\left( \widetilde{K_{m}}\right)
_{z}\right) , \left( \widetilde{K_{m}}\right) _{z}}. \label{BerTop}
\end{align}
By expressing \eqref{BerTop} as an integral,
% the Berezin transformation associated with the Hilbert space $ \gBerg $  is the integral operator
we get
\begin{align}\label{Berezin1}
\Berm \left[ \phi \right] (z) =\int_{\B}B_{m}(z,w) \phi (w) (1-|w|^2)^{-n-1}d\mu (w),
\end{align}
where
\begin{align}
B_{m}(z,w) :=\frac{\left| K_{m}^{\chp,n}(z,w) \right|^2%
}{K_{m}^{\chp,n}(z,z) }, \qquad z,w\in \B.
\end{align}
Explicitly, %since $P_{m}^{( n-1,2[\chp-m] -n) }(1)= \frac{\Gamma(n+m)}{m!\Gamma(n)}$, we infer
\begin{align}
B_{m}(z,w) = \frac{m!\Gamma(n)}{\Gamma(m+n)}\gamma^{\chp,n}_m (\cosh d(z,w))^{-4(\chp-m)}\left( P_{m}^{\left( n-1,2[\chp-m] \right) }\left( 1-2\tanh^2(d(z,w))
\right) \right)^2,
\end{align}
where $\gamma^{\chp,n}_m$ is as in \eqref{crkb}.
%where $\cosh^2(d(z,w))$ is as in \eqref{Dist}.
We then state the following

\begin{definition} Let $m=0,1,2, \cdots, [\chp-\frac n2]$ and let $\gamma^{\chp,n}_m$ as in \eqref{crkb}. Then, the integral operator acting on $\phi\in L^{\infty }\left( \B%
\right) $  as
\begin{align}
\Berm [\phi] (z) =&\frac{m!\Gamma(n)}{\Gamma(m+n)}\gamma^{\chp,n}_m    \int_{\B}(\cosh(d(z,w)))^{-4(\chp-m)}\label{GBT}\\
&\times
\left( P_{m}^{\left( n-1,2[\chp-m] \right) }\left( 1-2\tanh^2(d(z,w))
\right) \right)^2 \phi(w) \left( 1-|w|^2\right)^{-n-1}d\mu(w).\nonumber
\end{align}
 is called the generalized Berezin transform of index $m$.
\end{definition}

\begin{proposition} The transform $\Berm $ defined in \eqref{GBT} can be written as
\begin{align}
\Berm [\phi](z) &= \frac{m!\Gamma(n)}{\Gamma(m+n)}\gamma^{\chp,n}_m \int_{\B}(1-|w|^2)^{2(\chp-m)}  \label{RHS}\\
&\times \left( P_{m}^{( n-1,2[\chp-m])}  ( 1-2|w|^2) \right)^2 \phi\left( \varphi _{z}(w) \right) (1-|w|^2)^{-n-1} d\mu(w), \nonumber
\end{align}
where $\varphi _{z}(w)$ is as in \eqref{varphiz}.
\end{proposition}

\noindent {\it Proof.} Fix $z\in \B$ and set $w=\varphi_{z}^{-1}(\xi) $, where $\varphi _{z}(w)$ is as in \eqref{varphiz}. We have
$\varphi _{z}\in Aut(\B)$ and $\varphi _{z}\left( 0\right) =z$. The quantity in the right hand side of \eqref{RHS} reads, in terms of $\xi $, as
\begin{eqnarray}
&& \frac{m!\Gamma(n)}{\Gamma(m+n)}\gamma^{\chp,n}_m \left( 1-\left| \varphi _{z}^{-1}(\xi) \right|^2\right)^{2( \chp-m) }\label{expression} \\
&& \qquad \times
\left( P_{m}^{(n-1,2[\chp-m])}\left( 1-2\left|\varphi_{z}^{-1}(\xi)\right|^2\right) \right)^2 \phi(\xi)
\left(1-\left|\varphi_{z}^{-1}(\xi)\right|^2\right)^{-n-1} d\mu \left( \varphi _{z}^{-1}(\xi) \right). \nonumber
\end{eqnarray}
By writing down the identity \eqref{IdentityInv} for $\psi=\varphi_{z}^{-1}$ and $w=w'$, we get %
 \begin{equation}
1-\left| \varphi _{z}^{-1}(\xi) \right|^2= \frac{(1-|z|^2)
\left( 1-|\xi|^2\right) }{\left| 1-\scal{\xi,z} \right|^2} =\cosh^{-2}(d(z,\xi)) = 1-\tanh^2(d(z,\xi))
. \end{equation}
Also, by making use of \eqref{Jacobian} giving the explicit expression of
the Jacobian of the transformation $\xi =\varphi _{z}(w) $ is given through \eqref{Jacobian}, we see that
\begin{align}
\left( 1-\left| \varphi _{z}^{-1}(\xi) \right|^2\right) ^{-n-1} d\mu \left( \varphi _{z}^{-1}(\xi) \right) &= \left( 1-\left| \varphi _{z}^{-1}(\xi) \right|
^2\right) ^{-n-1} \frac{\left|Dw\right| }{\left| D\xi \right| }d\mu (\xi)\nonumber \\
&= (1-|\xi|^2)^{-n-1}d\mu (\xi). \label{Jacobian}
\end{align}
Therefore, the expression in \eqref{expression} can also be written as
\begin{eqnarray}
&& \frac{m!\Gamma(n)}{\Gamma(m+n)}\gamma^{\chp,n}_m \int_{\B} \left( \frac{(1-|z|^2) \left( 1-|\xi|^2\right) }{\left| 1-\scal{\xi ,z} \right|
^2}\right) ^{2(\chp-m)} \\ && \qquad \qquad\times \left( P_{m}^{\left( n-1,2[\chp-m] \right) }\left(
1-2\tanh^2(d(z,\xi)) \right) \right)^2 \phi(\xi) \left( 1-|\xi|^2\right) ^{-n-1}d\mu
(\xi) \nonumber
\end{eqnarray}
or equivalently as
\begin{align}
&\frac{m!\Gamma(n)}{\Gamma(m+n)}\gamma^{\chp,n}_m \int_{\B%
}\left( \cosh d(z,\xi) \right) ^{-4\left( \chp-m\right) }
\\ & \qquad \times \left(P_{m}^{\left( n-1,2[\chp-m] \right) }\left( 1-2\tanh^2(d(z,\xi))\right) \right)^2
\phi(\xi) \left( 1-|\xi|^2\right)^{-n-1}d\mu (\xi), \nonumber
\end{align}
which is exactly the expression of $\Berm $ as defined in \eqref{GBT}. \fin

\section{$\Berm $ as function of the Laplace-Beltrami operator}

Our task here is to express the generalized Berezin transform $\Berm $
 as a function of the Laplace-Beltrami operator $\Delta_{\B}$. Precisely, we establish the following result.
 % look for a formula of the function $f$ such that $\Berm  =  f(-\Delta_{\B} -n^2)$.  Namely, we have
\begin{theorem}
Let $m=0,1,2, \cdots, [\chp-\frac n2]$. Then, the transform $\Berm $ defined in \eqref{GBT}
%and set $\Lambda=\frac 12\sqrt{-\Delta_{\B} -n^2}$,
can be expressed in terms of the Laplace-Beltarmi operator $\Delta_{\B}$ as
\begin{align*}
\Berm  &=C^{\chp,n}_m \sum\limits_{j=0}^{2m}   (-2)^{j} A_{j}
    \frac{\Gamma(n+j)\Gamma\left(2(\chp-m)-\frac 12\left(n-i\sqrt{-\Delta_{\B} -n^2}\right) \right) }
    {\Gamma\left( 2(\chp-m) +j + \frac 12\left(n+i\sqrt{-\Delta_{\B} -n^2}\right)\right) }
    %\mathcal{B}\left( n+j,-\frac{n}{2}+\frac i2\sqrt{-\Delta_{\B} -n^2}+2(\chp-m) \right)
\\ &  \times  {\bF }\left[ \begin{array}{c}
                 \frac 12\left(n+i\sqrt{-\Delta_{\B} -n^2}\right), \, n+j, \, \frac 12\left(n+i\sqrt{-\Delta_{\B} -n^2}\right) \\
                 (\chp-m) +j+\frac 12\left(n+i\sqrt{-\Delta_{\B} -n^2}\right), \, n
               \end{array} \bigg| 1
\right], \nonumber
\end{align*}
where
\begin{align*}
C^{\chp,n}_m&= \frac{(2(\chp-m)-n) \Gamma(n+m)}{m!\Gamma(2\chp -n-m+1) \Gamma(2\chp -n)} \\
A_{j}&=2^{-j}\sum_{p=\max \left( 0,j-m\right) }^{\min \left( m,j\right) }%
%\binom{m}{j-p}\binom{m}{p}
\frac{(m!)^2\Gamma(2\chp-m) \Gamma(2\chp-m+j-p)}{(j-p)!(m+p-j)!p!(m-p)!\Gamma(n+j-p) \Gamma(n+p) }.
\end{align*}
\end{theorem}

\noindent{\bf Proof.}
On one hand, the transform $\Berm $ is the integral operator
\begin{align}\label{GBT2}
\Berm \left[ \varphi \right] (z) =\int_{\B}B_{m}(z,w) \varphi (w)  (1-|w|^2)^{-n-1}d\mu (w)
\end{align}
with the kernel function
%$B_{m}(z,w)$ depends only on the geodesic distance $\rho:=d(z,w)$ and is given explicitly through \eqref{GBK} as
\begin{align}\label{KerGBT}
%B_{m}(z,w) = \frac{\pi^{n}\Gamma \left( 2\chp -n-m+1\right) }{\left( 2\chp -n-2m\right) \Gamma \left( 2\chp -m\right) }  \frac{m!}{(n) _{m}} (\cosh d(z,w))^{-4(\chp-m)} \left(P_{m}^{\left( n-1,2\left[ \chp-m\right] \right) }\left( 1-2\left|w \right|^2\right) \right)^2.
B_{m}(z,w) = \frac{{m!}\Gamma(n)}{\Gamma(n+m)} \gamma^{\chp,n}_m
(\cosh (\rho))^{-4(\chp-m)} \left(P_{m}^{\left( n-1,2\left[ \chp-m\right] \right) }\left( 1-2\tanh^2 (\rho)\right) \right)^2,
\end{align}
where $\rho=d(z,w)$ and $\gamma^{\chp,n}_m$ is as in \eqref{crkb}.
%\begin{eqnarray}
%\gamma^{\chp,n}_m= \frac{(2(\chp-m) -n)\Gamma(2\chp -m)}{\pi^n\Gamma(2\chp -m-n+1)}.\label{crkb}
%\end{eqnarray}
On the other hand, from \eqref{FctLn} combined with \eqref{SpecKernel}, we obtain
\begin{align}
\label{FctLn2}
f(-\Delta_{\B} -n^2)[\varphi](z) = \int_{\B} \left(\int_{0}^{+\infty} \Psi(z,w;\lambda) f(\lambda) d\lambda \right) \varphi(w)  (1-|w|^2)^{-n-1}d\mu(w),
\end{align}
$\Psi(z,w;\lambda)$ being the spectral function given in \eqref{BLBOueld}.
%\begin{equation}
%\Psi(z,w;\lambda) =   \frac{|\Gamma(\frac{n+i\lambda}2)|^4}{4\pi^{n+1}\Gamma(n)|\Gamma(i\lambda)|^2}
%\ \ {\aF}\left( \frac{n+i\lambda }{2},\frac{n-i\lambda }{2}; n; -\sinh^2 (\rho)\right) \label{BLBOueld2}
%\end{equation}
%is the kernel function given in \eqref{BLBOueld}.
Now, by equating \eqref{GBT2} and \eqref{FctLn2} we get
\begin{equation}
\int_{0}^{+\infty }\Psi(z,w;\lambda) f(\lambda) d\lambda = B_{m}(z,w)   \label{BLB10}
\end{equation}
by uniqueness of the kernels.
Explicitly,
\begin{align}
\int_{0}^{+\infty }
& \frac{\left|\Gamma\left(\frac{n+i\lambda}{2}\right)\right|^4}{\left| \Gamma(i\lambda) \right|^2}\ \
 \aF \left[ \begin{array}{c} \frac{n+i\lambda }{2},\frac{n-i\lambda }{2}\\ n \end{array} \bigg| -\sinh^2(\rho) \right] f(\lambda) d\lambda =  h(\rho),  %\label{BLB14}
\label{BLB15}
\end{align}
 where $h$ is used to represents the function $4\pi^{n+1}{m!}\Gamma(n) B_{m}(z,w)$, i.e.,
 \begin{align}
  h(\rho):= %4 \frac{\pi^{2n+1}m!\Gamma^2(n)\Gamma(2\chp -n-m+1)}{(2\chp -n-2m) \Gamma(n+m)\Gamma(2\chp -m)}
 \frac{4\pi^{n+1}{m!}\Gamma^2(n)}{ \Gamma(n+m)}{\gamma^{\chp,n}_m}\left( \cosh (\rho) \right) ^{-4(\chp-m) } \left( P_{m}^{\left( n-1,2(\chp-m) -n\right)}(1-2\tanh^2(\rho)) \right)^2.  \label{h}
\end{align}
 To  determine the function $f$ we need to invert the integral equation \eqref{BLB15}. For this, we make appeal to the Fourier-Jacobi transform
 \begin{align}
 h\in L^2\left( \R_{+},\Delta _{\alpha ,\beta }(t) dt\right)  \longmapsto  g\in L^2\left( \R_{+},(2\pi) ^{-1}\left|c_{\alpha ,\beta }(\lambda) \right| ^{-2}d\lambda \right)
\end{align}
 as defined in \eqref{JacTransf} by
\begin{equation}
g(\lambda) =\int_{0}^{+\infty }h(t) \phi_{\lambda }^{(\alpha,\beta) }(t) \Delta _{\alpha
.\beta }(t) dt. \label{BLB17}
\end{equation}
Its inverse is given by
\begin{equation}
h(t) =\frac{1}{2\pi }\int_{0}^{+\infty }g(\lambda) \phi _{\lambda }^{(\alpha,\beta) }(t)
\left| c_{\alpha ,\beta }(\lambda) \right| ^{-2}d\lambda .
\label{BLB19}
\end{equation}
For $\alpha =n-1$, $\beta =0$ and  $t=\rho $, the involved quantities
  $\phi_{\lambda}^{(\alpha,\beta)}(t)$,
 $\Delta_{\alpha,\beta}(t)$ and $c_{\alpha ,\beta}(\lambda)$ given respectively by \eqref{JacFun}, \eqref{Deltaab}
 and \eqref{ChandraT} read
 \begin{align}
&\phi _{\lambda }^{(n-1,0) }(\rho) =\aF \left[ \begin{array}{c}
\frac{n+i\lambda}{2} ,\frac{1}{2}(n-i\lambda)
\\ n \end{array} \bigg | -\sinh^2(\rho) \right]   \label{C1}
\\
& \Delta_{n-1,0}(\rho) = 2 \cosh(\rho)(\sinh(\rho))^{2n-1} \label{C2}\\
&c_{n-1,0}(\lambda) =\frac{2^{n-i\lambda }\Gamma(n)
\Gamma(i\lambda) }{\Gamma^2\left( \frac{n+i\lambda}{2} \right)}.  \label{C3}
\end{align}
%so that we have
%\begin{equation}
%\left| c_{n-1,0}(\lambda) \right| ^{-2}=\frac{1}{2^{2n}\Gamma^2(n)}\frac{\left| \Gamma \left( \frac{n+i\lambda}{2} \right) \right| ^{4}}{\left| \Gamma(i\lambda) \right|
%^2}.  \label{BLB25}
%\end{equation}
Therefore, \eqref{BLB15}  can be rewritten in terms of $\phi_{\lambda}^{(n-1,0)}(\rho)$ and
 $c_{n-1,0}(\lambda)$  as
%\begin{equation}
%\int_{0}^{+\infty }2^{2n} \Gamma^2(n) \left| c_{n-1,0}(\lambda) \right| ^{-2}\phi _{\lambda
%}^{(n-1,0) }(\rho) f(\lambda) d\lambda
%=h(\rho)   \label{BLB27}
%\end{equation}
%or equivalently
\begin{equation}
\frac{1}{2\pi }\int_{0}^{+\infty }f(\lambda) \phi _{\lambda
}^{(n-1,0) }(\rho) \left| c_{n-1,0}(\lambda) \right|^{-2}d\lambda
=\frac{1}{2\pi  2^{2n}\Gamma^2(n) } h(\rho) . \label{BLB28}
\end{equation}
In view of the transform \eqref{BLB17}, the function $f$ takes the form
\begin{align}
f(\lambda) =\frac{1}{2 \pi 2^{2n}  \Gamma^2(n) }
\int_{0}^{+\infty} h(\rho) \phi_{\lambda }^{(n-1,0)}(\rho) \Delta_{n-1,0}(\rho) d\rho   . \label{BLB33}
\end{align}
 Substitution of $h(\rho)$, $\phi_{\lambda }^{(n-1,0)}(\rho)$ and
 $\Delta_{n-1,0}(\rho)$ by their explicit expressions, as given by \eqref{h}, \eqref{C1} ad \eqref{C2} respectively, yields
\begin{align}
f(\lambda) = & %\frac{2m!\left( 2\left( \chp-m\right) -n\right) \Gamma \left(
%2\chp -m\right) }{\Gamma \left( m+n\right) \Gamma \left( 2\chp -\left( n+m\right)+1\right) }
\frac{2\pi^{n}m!}{\Gamma(n+m)} \gamma^{\chp,n}_m\int_{0}^{+\infty }
{\aF}\left[ \begin{array}{c}  \frac{n+i\lambda}{2},\frac{n-i\lambda}{2}\\n \end{array} \bigg |-\sinh^2(\rho) \right] \label{flamb}
\\ & \times  (\sinh(\rho))^{2n-1}\left( \cosh(\rho) \right) ^{-4(\chp-m)
+1} \left( P_{m}^{(n-1,2(\chp-m)-n)}(1-2\tanh^2(\rho)) \right)^2d\rho . \nonumber
\end{align}
By the change of variable $x=\sinh^2(\rho)$, Equation \eqref{flamb} becomes
\begin{align}
f(\lambda)% &=C_{\chp,m,n}\int_{0}^{+\infty }F_{1}\left(
%\frac{n+i\lambda }{2},\frac{n-i\lambda }{2},n,-x\right) \left( x^{\frac{1}{2}%
%}\right) ^{2n-1}\left( 1+x\right) ^{-2(\chp-m) }\left( 1+x\right) ^{%
%\frac{1}{2}} \nonumber \\ & \qquad\qquad
%\times \left( P_{m}^{\left( n-1,2(\chp-m) -n\right) }\left( \frac{%
%1-x}{1+x}\right) \right)^2\frac{1}{2}x^{-\frac{1}{2}}\left( 1+x\right) ^{-%
%\frac{1}{2}}dx \nonumber \\
%&=\frac{1}{2}\frac{2m!\left( 2\left( \chp-m\right) -n\right) \Gamma \left(
%2\chp -m\right) }{\Gamma \left( m+n\right) \Gamma \left( 2\chp -\left( n+m\right)+1\right) }
=\frac{\pi^{n}m!}{ \Gamma(n+m)} \gamma^{\chp,n}_m& \int_{0}^{+\infty }{\aF} \left[ \begin{array}{c} \frac{%
n+i\lambda }{2},\frac{n-i\lambda }{2}\\n \end{array} \bigg | -x\right] \label{BLBBack}
\\  & \times x^{n-1}(1+x)^{-2(\chp-m)} \left( P_{m}^{\left( n-1,2(\chp-m) -n\right)
}\left( \frac{1-x}{1+x}\right) \right)^2dx .\nonumber
\end{align}
Next, we adopt a formula of the product of two Jacobi polynomials \cite[Eq (10), p.648]{TT} to write that
 \begin{equation}
\left( P_{m}^{(n-1,2(\chp-m)-n) }\left( \frac{1-x}{1+x}\right) \right)^2
=\frac{\left( \Gamma(n+m) \right)^2}{(m!)^2\left( \Gamma(2\chp -m) \right)^2}
\sum\limits_{j=0}^{2m}(-1) ^{j}A_{j}\left( 1-\frac{1-x}{1+x}\right) ^{j}, \label{JacProd}
 \end{equation}
where
 \begin{equation}\label{Aj}
A_{j}=2^{-j}\sum_{p=\max \left( 0,j-m\right) }^{\min \left( m,j\right) }%
%\binom{m}{j-p}\binom{m}{p}
\frac{(m!)^2\Gamma(2\chp-m) \Gamma(2\chp-m+j-p)}
{(j-p)!(m+p-j)!p!(m-p)!\Gamma(n+j-p) \Gamma(n+p) }.
\end{equation}
Inserting \eqref{JacProd} into \eqref{BLBBack}, we arrive at
\begin{align}
f(\lambda) =%\frac{\left( 2\left( \chp-m\right) -n\right) \Gamma( n+m)}{\Gamma( 2\chp -(n+m) +1) \Gamma( 2\chp -n) }
 \pi^{n}\gamma^{\chp,n}_m \Gamma(n+m)  \sum\limits_{j=0}^{2m}\left( -2\right) ^{j}A_{j} I_{n,j}(\lambda),
\end{align}
where
\begin{align}
I_{n,j}(\lambda) :=\int_{0}^{+\infty }\ \
\aF\left[ \begin{array}{c} \frac{n+i\lambda }{2},\frac{n-i\lambda }{2}\\ n \end{array} \bigg | -x\right]
x^{n-1+j}\left( 1+x\right) ^{-2(\chp-m) -j}dx.
\end{align}
Using the integral representation for the $\aF$-sum \cite[p.1005]{GR}:
\begin{align}
{\aF}\left[ \begin{array}{c} \alpha ,\beta \\ \gamma \end{array} \Big | z\right]  =\frac{1}{\mathcal{B}\left( \beta ,\gamma -\beta \right) }
\int_{0}^{1}t^{\beta -1}(1-t)^{\gamma -\beta -1}
                                      (1-tz)^{-\alpha }dt, \quad
                                      \Re\left( \gamma \right) >\Re \left( \beta \right) >0 ,
\end{align}
where $\mathcal{B}(\cdot,\cdot)$ is the Beta function, for the special case of
 $ \alpha =\frac{n+i\lambda }{2}$, $\beta =\frac{n-i\lambda }{2}$, $\gamma =n$ and $z=-x$, it follows
\begin{align}
I_{n,j}(\lambda) & = \int_{0}^{\infty }\left[ \frac{\Gamma(n) }{\left| \Gamma
\left( \frac{n+i\lambda }{2}\right) \right|^2}\int_{0}^{1}t^{\frac{%
n-i\lambda }{2}-1}(1-t) ^{\frac{n+i\lambda }{2}-1}\left(
1+tx\right) ^{-\frac{\left( n+i\lambda \right) }{2}}dt\right]
x^{n-1+j}\left( 1+x\right) ^{-2\left( \chp-m\right) -j}dx \nonumber\\
&=\frac{\Gamma(n) }{|\Gamma(\frac{n+i\lambda }{2})|^2}
\int_{0}^{1} t^{\frac{n-i\lambda }{2}-1}(1-t)^{\frac{n+i\lambda }{2}-1}
Q_{\chp,m,j,\lambda }^{n}(t)  dt , \label{IntI1}
\end{align}
where
\begin{align}
Q_{\chp,m,j,\lambda }^{n}(t) = \int_{0}^{\infty }\left( 1+tx\right)^{-\frac{( n+i\lambda) }{2}}
              x^{n-1+j}(1+x)^{-2(\chp-m) -j}dx . \label{IntQj0}
\end{align}
Making use of the identity \cite[p.317]{GR}:
\begin{align}
\int_{0}^{\infty }x^{\lambda_{0}-1}(1+x)^a(1+\alpha x)^{b}dx =
\mathcal{B}(\lambda_{0},-(a+b)-\lambda_{0}) \aF \left[ \begin{array}{c} -b ,\lambda_{0}\\ -(a+b) \end{array} \Big | 1-\alpha \right]
\end{align}
requiring that $|\arg(\alpha)| <\pi$ and $-\Re(a+b) >\Re(\lambda_{0})>0$ that are both fulfilled when $\lambda _{0}=n+j$,
$a =-2\left( \chp-m\right) -j$, $b=-(n+i\lambda)/2$ and $\alpha =t$, we get
\begin{align}
Q_{\chp,m,j,\lambda }^{n}(t) =\mathcal{B}\left( n+j,\frac{i\lambda -n}{2}+2\left( \chp-m\right) \right) \ \
\aF\left[ \begin{array}{c} \frac{n+i\lambda }{2},n+j\\ \frac{n+i\lambda }{2}+2\left( \chp-m\right) +j\end{array} \bigg| 1-t\right]. \label{IntQj}
\end{align}
By inserting \eqref{IntQj} in  \eqref{IntI1}, we see that
\begin{align}
I_{n,j}(\lambda) &=\frac{\Gamma(n)}{\left| \Gamma\left(\frac{n+i\lambda}{2}\right) \right|^2}
  \mathcal{B}\left( n+j,
\frac{i\lambda -n}{2}+2(\chp-m) \right) \\ &\times \int_{0}^{1}t^{\frac{%
n-i\lambda }{2}-1}(1-t) ^{\frac{n+i\lambda }{2}-1} {\aF}
\left[ \begin{array}{c} \frac{n+i\lambda }{2},n+j\\ \frac{n+i\lambda }{2} +2(\chp-m) +j\end{array} \bigg |1-t\right] dt.\nonumber
\end{align}
The change of the variable $1-t=u$ and the use of the formula \cite[p.813]{GR}:
 \begin{equation}
\int_{0}^{1}u^{a-1}(1-u)^{b-1}\aF \left[ \begin{array}{c} \alpha ,\beta \\ \gamma \end{array} \Big | u\right] du =\frac{\Gamma(a)\Gamma(b)}{\Gamma(a+b)}\ \
{\bF }\left[ \begin{array}{c} \alpha ,\beta ,a \\ \gamma , a+b \end{array} \Big| 1 \right],
 \end{equation}
$\Re(a) >0$, $\Re(b) >0$ and $\Re(b+\gamma -\alpha -\beta) >0$ for
$a = \frac{n+i\lambda }{2}$, $b =\frac{n-i\lambda }{2}$, $\alpha =\frac{n+i\lambda
}{2},\beta =n+j$ and $\gamma =\frac{n+i\lambda }{2}+2\left( \chp-m\right) +j$, yield
\begin{align}
I_{n,j}(\lambda) &=B\left( n+j,-\frac{n}{2}+\frac{i\lambda }{2}+2(\chp-m) \right) \label{Ij}
\\ & \qquad \times {\bF }\left[ \begin{array}{c} \frac{n+i\lambda }{2},n+j,\frac{n+i\lambda }{2}\\
\frac{n+i\lambda }{2}+2(\chp-m) +j,n\end{array} \Big| 1 \right]. \nonumber
\end{align}
 Finally, we arrive at
\begin{align}
f(\lambda) &= %d_{\chp,m,n}
\frac{(2(\chp-m) -n) \Gamma(n+m)
}{m!\Gamma(2\chp-(n+m)+1) \Gamma(2\chp-m) }\sum\limits_{j=0}^{2m}(-2)^{j} \mathcal{B}\left( n+j,\frac{i\lambda-n}{2}+2(\chp-m) \right)
A_{j} \nonumber
\\ & \qquad \times {\bF }\left[ \begin{array}{c} \frac{n+i\lambda }{2},n+j,\frac{
n+i\lambda }{2}\\ \frac{n+i\lambda }{2}+2(\chp-m) +j,n \end{array} \Big | 1\right].\label{fm}
\end{align}
Replacing $\lambda$ by $\sqrt{-\Delta_{\B}-n^2}$ and expressing the Beta function in terms of Gamma functions we obtain the announced formula. \fin

\section{ The case $m=0$}

In this case, the transform $\Berezin  $ turns out to be the well known Berezin transform $\mathscr B$
related to the Bergman ball of $\Bbb{C}^{n}$
\begin{equation}
\mathscr B \left[ \varphi \right] (z) = \frac{(2\chp -n)\Gamma(2\chp )}{\pi^n\Gamma(2\chp -n+1)}
\int_{\B}\left(\frac{(1-|z|^2)(1-|\xi|^2)}{|1-<z,\xi>|^2}\right)^{2\chp } \varphi
(\xi) \left( 1-|\xi|^2\right) ^{-\left( n+1\right)
}d\mu(\xi)   \label{BLB3}
\end{equation}
with the kernel function
 \begin{equation}
B(z,w) = \frac{(2\chp -n)\Gamma(2\chp )}{\pi^n\Gamma(2\chp -n+1)}
\left(\frac{(1-|z|^2)(1-|\xi|^2)}{|1-<z,\xi>|^2}\right)^{2\chp }.
 \end{equation}
Based on a formula due to Berezin \cite[p.377]{Berezin2}, a formula  expressing $\mathscr B  $ as a function of the Laplace-Beltrami
operator $\Delta_{\B}$ has been derived by  Peeter \cite[p.182]{peeter} as
 \begin{equation}\label{Lambdaa}
\mathscr B  =\frac{\left| \Gamma\left(\alpha+1+\frac{n}{2} +\frac i2 \sqrt{-\Delta_{\B}-n^2} \right)
\right|^2}{\Gamma(\alpha+1) \Gamma \left( \alpha +n+1\right) }.
 \end{equation}
 There $\alpha =\frac{1}{\hbar }-n-1$ occurs in the measure weight of the Bergman
space considered in \cite{peeter}. The result we obtain recovers the Peeter's one by taking $m=0$ and $\alpha=2\chp-n-1$ (i.e., $\frac{1}{\hbar }=2\chp$). Indeed, we establish the following

\begin{proposition}
 For $m=0$, the function in \eqref{fm} reduces to
 \begin{equation}
f_{0}(\lambda) =\frac{\left| \Gamma(2\chp -\frac{n-i\lambda}{2}) \right|^2}{\Gamma(2\chp-n) \Gamma(2\chp) }.
 \end{equation}
\end{proposition}

\noindent {\bf Proof.}
For $m=0$, the expression in \eqref{fm} reads
\begin{align}
f_0(\lambda) &= \frac{(2\chp-n)\Gamma(n)}{\Gamma(2\chp-n+1)\Gamma(2\chp)}
                \frac{\Gamma(n)\Gamma\left(\frac{i\lambda-n}{2}+2\chp \right)}{\Gamma\left(2\chp+\frac{n+i\lambda}{2} \right)}
                A_{0} \\
                &\qquad  \times {\bF }\left[ \begin{array}{c} \frac{n+i\lambda }{2},n,\frac{n+i\lambda }{2}\\ \frac{n+i\lambda }{2}+2\chp,n \end{array} \Big | 1\right].\nonumber
\end{align}
Now, according to the expression of $A_0$ given through \eqref{Aj} as
\begin{equation}
A_0=\left(\frac{\Gamma(2\chp)}{\Gamma(n)}\right)^2
\end{equation}
 together with the equality between the hypergeometric series $\bF$ and $\aF$:
\begin{equation}
\bF\left[\begin{array}{c} a,b,c\\ d,e \end{array} \Big | x\right]=\aF\left[ \begin{array}{c} a,b\\ d \end{array} \Big |x \right]
\end{equation}
in the particular case when $c=e$, it follows that
\begin{align}
f_0(\lambda) = \frac{(2\chp-n)\Gamma(2\chp)}{\Gamma(2\chp-n+1)}
                \frac{\Gamma\left(2\chp +\frac{i\lambda-n}{2}\right)}{\Gamma\left(2\chp+\frac{n+i\lambda}{2} \right)}
                \, \,\aF\left[ \begin{array}{c} \frac{n+i\lambda }{2},\frac{n+i\lambda }{2}\\ \frac{n+i\lambda }{2}+2\chp
                 \end{array} \bigg | 1\right] .\nonumber
\end{align}
Next, by applying the Gauss theorem \cite[p.1008]{GR}
 \begin{equation}
\aF\left[ \begin{array}{c} \alpha ,\beta \\ \gamma \end{array} \Big | 1\right] =
\frac{\Gamma(\gamma)\Gamma(\gamma -\alpha -\beta)}
{\Gamma(\gamma-\alpha)\Gamma(\gamma -\beta) },
 \end{equation}
whenever $\Re(\gamma-\alpha-\beta) >0$, $\gamma \neq 0,-1,-2, \cdots$, we arrive at
\begin{align*}
f_0(\lambda)  = \frac{(2\chp-n)  |\Gamma(2\chp-\frac{n-i\lambda}{2})|^2 }{\Gamma(2\chp-n+1)\Gamma(2\chp)}
                .\nonumber
\end{align*}
Making use of the identity $\Gamma(x+1)=x\Gamma(x)$ leads to
\begin{align*}
f_{0}(\lambda) =\frac{\left| \Gamma \left( 2\chp -\frac{n-i\lambda}{2}\right)
\right|^2}{\Gamma \left( 2\chp -n\right) \Gamma \left( 2\chp \right) }.
\end{align*}
This ends the proof. \fin

As consequence, replacing $\lambda$ by $\sqrt{-\Delta_{\B}-n^2}$ in the result of the previous proposition yields get \eqref{Lambdaa}.

\end{document}